\documentclass[a4paper, 10pt, conference]{ieeeconf}      

\IEEEoverridecommandlockouts                              
\usepackage{graphics} 
\usepackage{amsmath} 
\usepackage{amssymb}  
\usepackage{epsfig}
\usepackage{graphicx}
\usepackage[hang]{caption2}
\usepackage[normalsize,it,hang]{subfigure}

\newtheorem{remark}{Remark}[section]

\title{\LARGE \bf
``Intelligent'' controllers \\ on cheap and small \\ programmable devices
}

\author{C\'edric Join$^{a,b,d}$, Fr\'ed\'eric Chaxel$^{b}$ and Michel Fliess$^{a,c}$
\thanks{$^{a}$ AL.I.E.N. (ALg\`{e}bre pour Identification et Estimation Num\'{e}riques), 24-30 rue Lionnois, BP 60120, 54003 Nancy, France. \newline
 {\tt \scriptsize \{cedric.join, michel.fliess\}@alien-sas.com}}
\thanks{$^{b}$ CRAN (CNRS, UMR 7039), Universit\'e de Lorraine, BP 239, 54506 Vand{\oe}uvre-l\`es-Nancy, France.
 \newline       {\tt \scriptsize \{cedric.join, frederic.chaxel\}@univ-lorraine.fr}}%
\thanks{$^{c}$ LIX (CNRS, UMR 7161), \'Ecole polytechnique, 91128 Palaiseau, France.
 {\tt \scriptsize Michel.Fliess@polytechnique.edu}}%
 \thanks{$^{d}$ Projet Non-A, INRIA Lille -- Nord-Europe, France.
       } }%

\begin{document}

\maketitle
\thispagestyle{empty}
\pagestyle{empty}

\begin{abstract}
It is shown that the ``intelligent'' controllers which are associated to the recently introduced model-free control synthesis may be easily implemented on cheap and small programmable devices. 
Several successful numerical experiments are presented with a special emphasis on fault tolerant control. 

\end{abstract}

\section{Introduction}
It is well know that  the overwhelming majority of industrial control applications is based on PID controllers (see, \textit{e.g.}, \cite{astrom,od}, and the references therein). 
Those controllers are often manufactured by numerous companies as 
microcontrollers on cheap and small programmable devices, like a \emph{Microchip}\footnote{Let us quote Wikipedia: \emph{Microchip Technology} is an American manufacturer of microcontroller, memory and analog semiconductors.} 
\emph{PIC} or a \emph{Freescale}\footnote{Let us quote again Wikipedia: \emph{Freescale Semiconductor} is an American company that produces and designs embedded hardware.} \emph{DSC} for instance This communication demonstrates that the recent \emph{model-free control} \cite{ijc13} and the corresponding \emph{intelligent PID controllers}  \cite{ijc13}
may also be implemented on such devices, 
\begin{itemize}
\item thanks to the low power computation cost which they require,
\item although they obey to mathematical principles, which are quite different from those of ``classic'' PIDs.
\end{itemize} 
When compared to a PID regulator, its intelligent counterpart contains one term more, which subsumes not only the unknown structure of the plant
but also the unknown disturbances. This term, which is estimated online thanks to recent parameter identification techniques (\cite{sira1,sira2}), is also found in the linear \emph{ultra-local} model
which 
\begin{itemize}
\item replaces the unknown global description of the plant,
\item is continuously updated online via the same calculations,
\item is of low order, which is most of the time equal to $1$.
\end{itemize}
\begin{remark}  
Let us emphasize moreover that the strange ubiquity of classic PIDs was mathematically explained 
for the first time, to the best of our knowledge, thanks to intelligent PIDs \cite{ijc13}. 
\end{remark} 
This aim is fully justified by the following facts:
\begin{itemize}
\item Intelligent PIDs are much easier to tune than the classic ones. 
\item They are robust with respect to most disturbances, including quite strong ones.
\item They permit straightforward fault accommodations.
\item Many successful concrete applications were already achieved in most various domains 
within a few years (see the numerous references in \cite{ijc13}). 
\end{itemize} 
Our paper is organized as follows. Section \ref{mfc} summarizes some of the most important facts about model-free control and its corresponding intelligent controllers. 
The implementation on small programmable devices is detailed in Section \ref{imple}. Section \ref{expe} describes some numerical experiments, with a peculiar emphasis on fault tolerant control, 
\textit{i.e.}, on an important topic in control engineering (see, {\it e.g.}, \cite{blanke,Noura}, and the references therein). Several excellent simulations are provided. Some 
concluding remarks may be found in Section \ref{conclu}.

\section{Model-free control: Basics\protect\footnote{See \cite{ijc13} for more details.}}\label{mfc}

\subsection{The ultra-local model}
The unknown global description of the plant is replaced by the \emph{ultra-local model}
\begin{equation}
{y^{(\nu)} = F + \alpha u} \label{ultralocal}
\end{equation}
where
\begin{itemize}
\item the derivation order $\nu \geq 1$ is selected by the practitioner;
\item $\alpha \in \mathbb{R}$ is chosen by the practitioner such that $\alpha u$ and
$y^{(\nu)}$ are of the same magnitude.
\end{itemize}
\begin{remark}
Note that $\nu$ has no connection with the order of the unknown system, which may be with distributed parameters, \textit{i.e.}, which might be best described by
partial differential equations (see, \textit{e.g.}, \cite{edf} for hydroelectric power plants).
\end{remark}
\begin{remark}
The existing
examples show that $\nu$ may always be chosen quite low,
\textit{i.e.}, $1$ or $2$. In almost all existing concrete case-studies $\nu = 1$. The only counterexample until now where $\nu = 2$ is provided by magnetic
bearings \cite{cifa12} where frictions are almost negligible.\footnote{See the explanation in \cite{ijc13}.}
\end{remark}
Some comments on $F$ are in order:
\begin{itemize}
\item $F$ is estimated via the measure of $u$ and $y$;
\item $F$ subsumes not only the unknown structure of the system but also
any perturbation.
\end{itemize}

\subsection{Intelligent PIDs}
Set $\nu = 2$ in Equation \eqref{ultralocal}:
\begin{equation}
\ddot{y} = F + \alpha u \label{2}
\end{equation}
Close the loop via the \emph{intelligent
proportional-integral-derivative controller}, or \emph{iPID},
\begin{equation}\label{ipid}
{u = - \frac{F - \ddot{y}^\ast + K_P e + K_I \int e + K_D
\dot{e}}{\alpha}}
\end{equation}
where 
\begin{itemize}
\item $e = y - y^\star$ is the tracking error,
\item $K_P$, $K_I$, $K_D$ are the usual tuning gains. 
\end{itemize}
Combining
Equations \eqref{2} and \eqref{ipid} yields
$$
\ddot{e} + K_D \dot{e} + K_P e + K_I \int e = 0
$$
where $F$ does not appear anymore. The tuning of $K_P$, $K_I$, $K_D$
is therefore quite straightforward. This is a major benefit when
compared to the tuning of ``classic'' PIDs (see, \textit{e.g.},
\cite{astrom,od}, and the references therein). 
\begin{remark}
If $K_I = 0$ we obtain the \emph{intelligent
proportional-derivative controller}, or \emph{iPD},
\begin{equation*}\label{ipd}
{u = - \frac{F - \ddot{y}^\ast + K_P e + K_D
\dot{e}}{\alpha}}
\end{equation*} 
\end{remark}

Set now $\nu = 1$ in Equation \eqref{ultralocal}:
\begin{equation}
\dot{y} = F + \alpha u \label{1}
\end{equation}
The loop is closed by \emph{intelligent proportional-integral
controller}, or \emph{iPI},
\begin{equation}\label{ipi}
u = - \frac{F - \dot{y}^\ast + K_P e + K_I \int e}{\alpha}
\end{equation}
If $K_I = 0$, it yields an \emph{intelligent
proportional controller}, or \emph{iP},
\begin{equation}\label{ip}
\boxed{ u = - \frac{F - \dot{y}^\ast + K_P e}{\alpha} }
\end{equation}
\begin{remark}
Equation \eqref{ip} and the corresponding iPs are most common in practice. This is again a major simplification 
with respect to ``classic'' PIDs and PIs.
\end{remark}

\subsection{Estimation of $F$}
$F$ in Equation \eqref{ultralocal} is assumed to be ``well'' approximated by a piecewise constant function $F_{\text{est}} $. According to the algebraic parameter identification 
developed in \cite{sira1,sira2}, rewrite, if $\nu = 1$, Equation \eqref{1}  in the operational domain (see, \textit{e.g.}, \cite{yosida}) 
$$s
Y = \frac{\Phi}{s}+\alpha U +y(0)
$$
where $\Phi$ is a constant. We get rid of the initial condition $y(0)$ by multiplying both sides on the left by $\frac{d}{ds}$:
$$
Y + s\frac{dY}{ds}=-\frac{\Phi}{s^2}+\alpha \frac{dU}{ds}
$$
Noise attenuation is achieved by multiplying both sides on the left by $s^{-2}$.\footnote{See \cite{bruit} for a theoretical explanation.} It yields in the time domain the realtime estimate
\begin{equation*}\label{integral}
{\small F_{\text{est}}(t)  =-\frac{6}{\tau^3}\int_{t-\tau}^t \left\lbrack (\tau -2\sigma)y(\sigma)+\alpha\sigma(\tau -\sigma)u(\sigma) \right\rbrack d\sigma }
\end{equation*}
where $\tau > 0$ might be quite small. This integral may of course be replaced in practice by a classic digital filter.


\section{Implementation}\label{imple}
Let us remind that implementing controllers on small programmable devices is a well established topic in engineering (see, \textit{e.g.}, \cite{godse,ibrahim,pa,tav,val,wil}, and the references therein). 
\subsection{The iP device}
We first detail the implementation of the iP device, which is most of the time
sufficient in practice.
\subsubsection{Device}\label{material}
The device is a microchip dsPIC33FJ128GP204 characterized by:\footnote{See {\tt http$:$//www.microchip.com} for more technical details, and also the prices.}
\begin{itemize}
\item Architecture:	16-bit.
\item CPU speed (MIPS):	40.
\item Memory type:	Flash.
\item Program memory (KB): 128.
\end{itemize}
We also utilize  
\begin{itemize}
\item two inputs  with a 12-bit, 500 KSPS analog-to-digital conversion,
\item one 12-bit digital output coupled with an external digital-to-anolog converter.
\end{itemize}
The values of the input and output variables, which are expressed in volts, are in the range $[0, 3.3]$. They can be connected however to an electronic amplifier stage. 
Let us add that our material is cheap: it costs less than $5$ euros.

Figure \ref{CS1} displays the corresponding architecture:
\begin{figure}[htb]
\begin{center}
{\resizebox*{7.5cm}{!}{\includegraphics{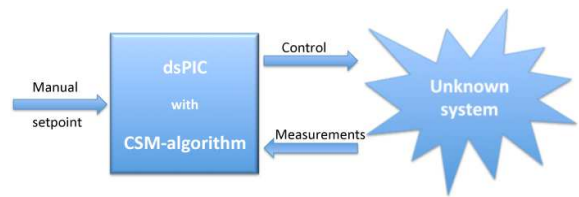}}}%
\caption{Control scheme}%
\label{CS1}
\end{center}
\end{figure}

\subsubsection{Simple calculations}
A control iteration needs:
\begin{itemize}
\item affectations: 19,
\item conditions (\textit{if...then...else}): 5,
\item additions/subtractions: 14,
\item multiplications/divisions: 16.
\end{itemize}
The computational power which is used is indeed very low.

\subsection{Some extensions}
Some new calculations are of course needed for iPIs, iPIDs, and iPDs.
\begin{remark}
An iPD was employed only once, for magnetic bearings \cite{cifa12}. Note moreover that 
it was not necessary until now to use iPIDs in practice!
\end{remark}

The integrals $\int e$ appearing in iPIs and iPIDs may be dealt via the classic trapezoidal rules (see, \textit{e.g.},  \cite{DI}). 
Less than $10$ basic operations are used.

The derivative $D(t)=K_D\dot e$ in iPIDs and iPDs may be obtained via
\begin{itemize}
\item a backward Euler difference scheme (see, \textit{e.g.},  \cite{DI})
$$D(t_k) = K_D \frac{e(t_k)-e(t_{k-1})}{t_k - t_{k-1}}$$
\item a low-pass digital filter for reducing the noise (see, \textit{e.g.}, \cite{DD}).
\end{itemize}
As above for the integration, less than 10 basic operations are used.
\begin{remark}
With very noisy signals, more advanced tools might be necessary for the differentiation (see \cite{nl,mboup}, and \cite{liu}).
\end{remark}

\section{Numerical experiments}\label{expe}

\subsection{LabVIEW}

\begin{figure}[htb]
\begin{center}
{\resizebox*{8.5cm}{!}{\includegraphics{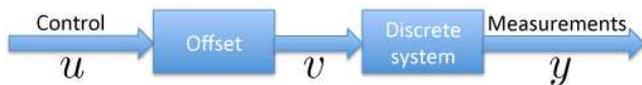}}}%
\caption{LabVIEW details}%
\label{Lab}
\end{center}
\end{figure}

As well known, LabVIEW greatly facilitates numerical simulations in engineering and in science.\footnote{See {\tt http$:$//www.ni.com/labview/whatis} for more details.} 
LabView can be programmed in order to control a real system by the use of an input-output device (analogic and logic inputs and outputs). Simple as well as advanced regulation strategies may therefore be assessed.
It may also be used for emulating/simulating the behavior of a real system (see, \textit{e.g.}, \cite{leger}). Emulation is achieved in our experimental platform 
via Labview and an input-output card (see Figure \ref{Lab}).

\begin{figure*}[htb]
\begin{center}
{\resizebox*{14.5cm}{!}{\includegraphics{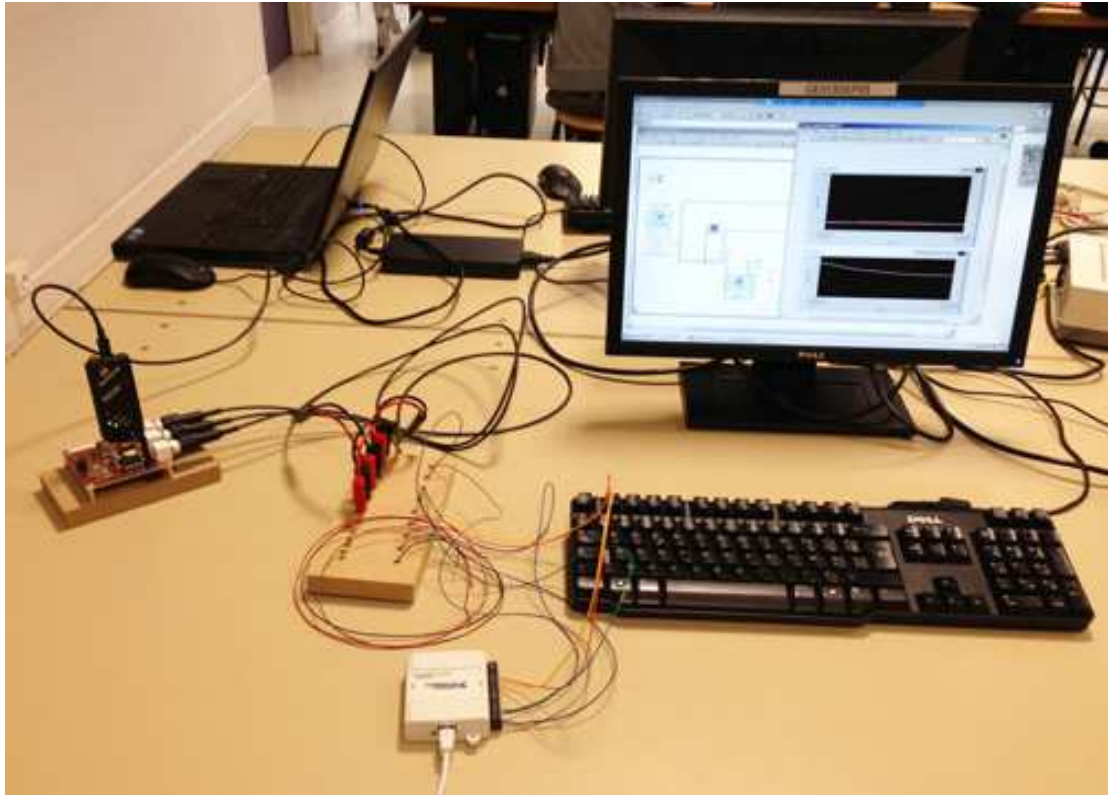}}}%
\caption{Our experimental setup with a dsPIC connected to an input-output card driven by LabVIEW}%
\label{plan}
\end{center}
\end{figure*}

Note that $v$ becomes the true input control variable, where $-1.65 \leq v \leq 1.65$.  The non-negativity condition on $u$ (see Section \ref{material}) is therefore dropped. Negative values for the control variables
are now possible.
Figure \ref{plan} displays a full description of our devices where the dsPIC is connected to an interface from National Instrument.\footnote{See {\tt http$:$//www.ni.com/white-paper/2732/en} for details.}

\subsection{Two experiments}
In both experiments the same iP is implemented where $\alpha=1$ and $K_P=1$. 
The control variable is obtained thanks to the dsPIC with a sample time equal to $0.1\textup{ms}$.

\subsubsection{System 1}
\begin{figure*}[htb]
\begin{center}
\subfigure[Control (-- blue) and control limits (- - red)]{
\resizebox*{8cm}{!}{\includegraphics{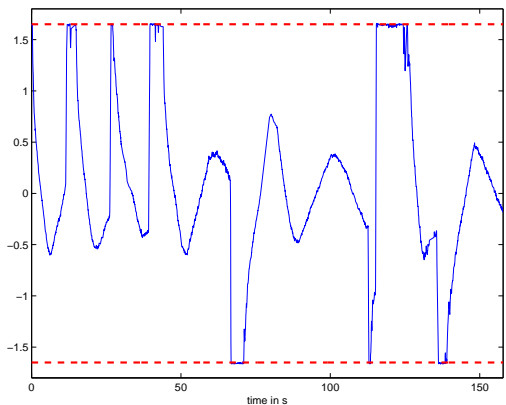}}}%
\subfigure[Setpoint (- - black) and output (-- blue)]{
\resizebox*{8cm}{!}{\includegraphics{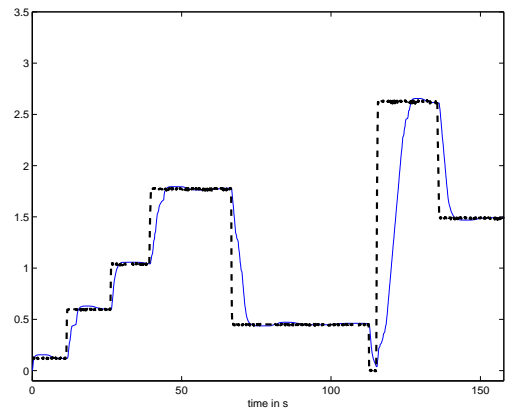}}}%
\caption{Experimental results: fault-free case}%
\label{rez1}
\end{center}
\end{figure*}

Consider the nonlinear stable system: 
$$
2(3\dot vv^2+v^3)=0.5\ddot y+0.5\dot y -y
$$
where $v$ is the saturated control after offset (see Figure \ref{Lab} for an explanation).
The tracking performances, which are presented in Figure \ref{rez1}, are quite good. They however deteriorate with important change points. This is due to the saturation of the control variable 
which is bounded.

\begin{figure*}[htb]
\begin{center}
\subfigure[Control (-- blue) and control limits (- - red)]{
\resizebox*{8cm}{!}{\includegraphics{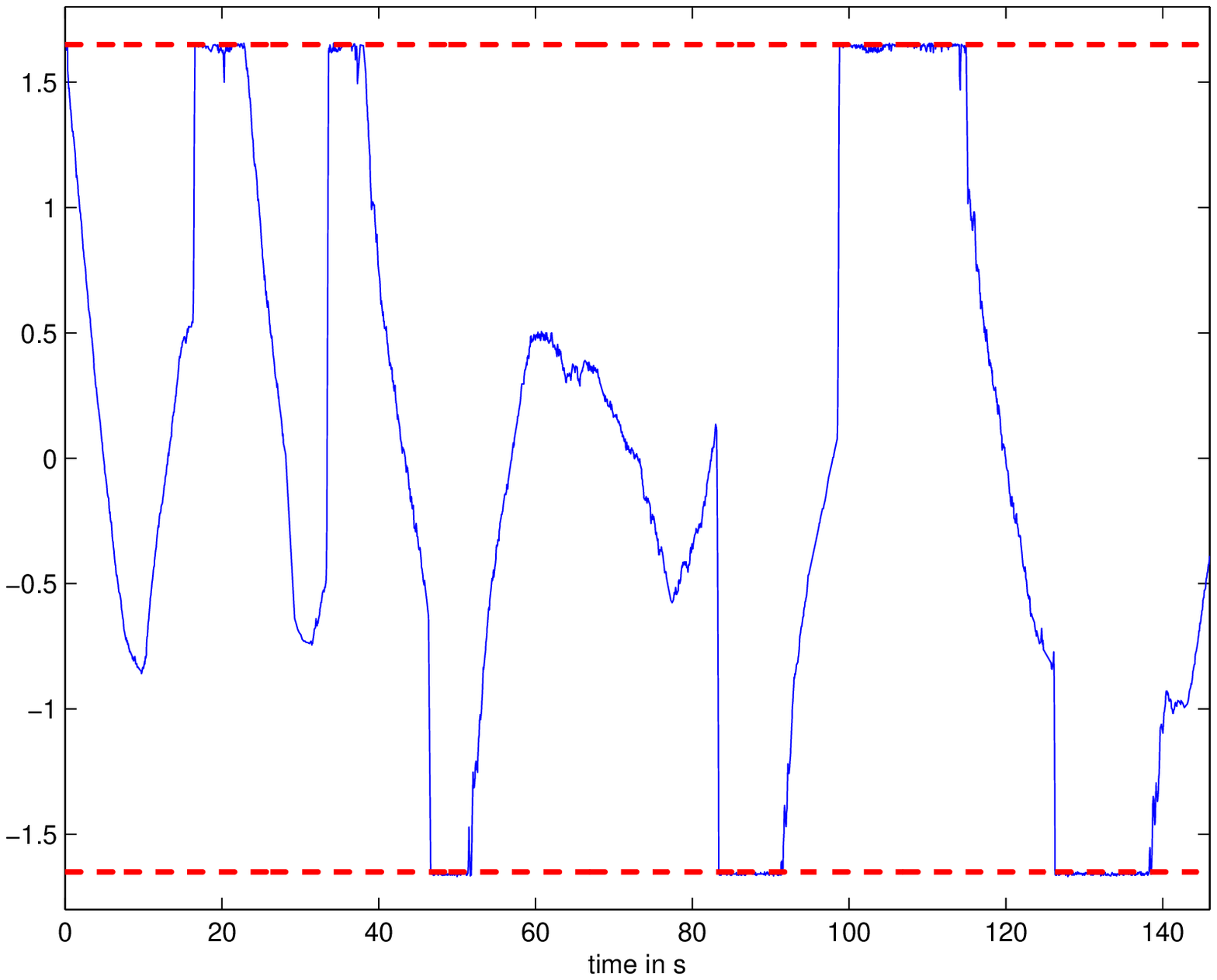}}}%
\subfigure[Setpoint (- - black) and output (-- blue)]{
\resizebox*{8cm}{!}{\includegraphics{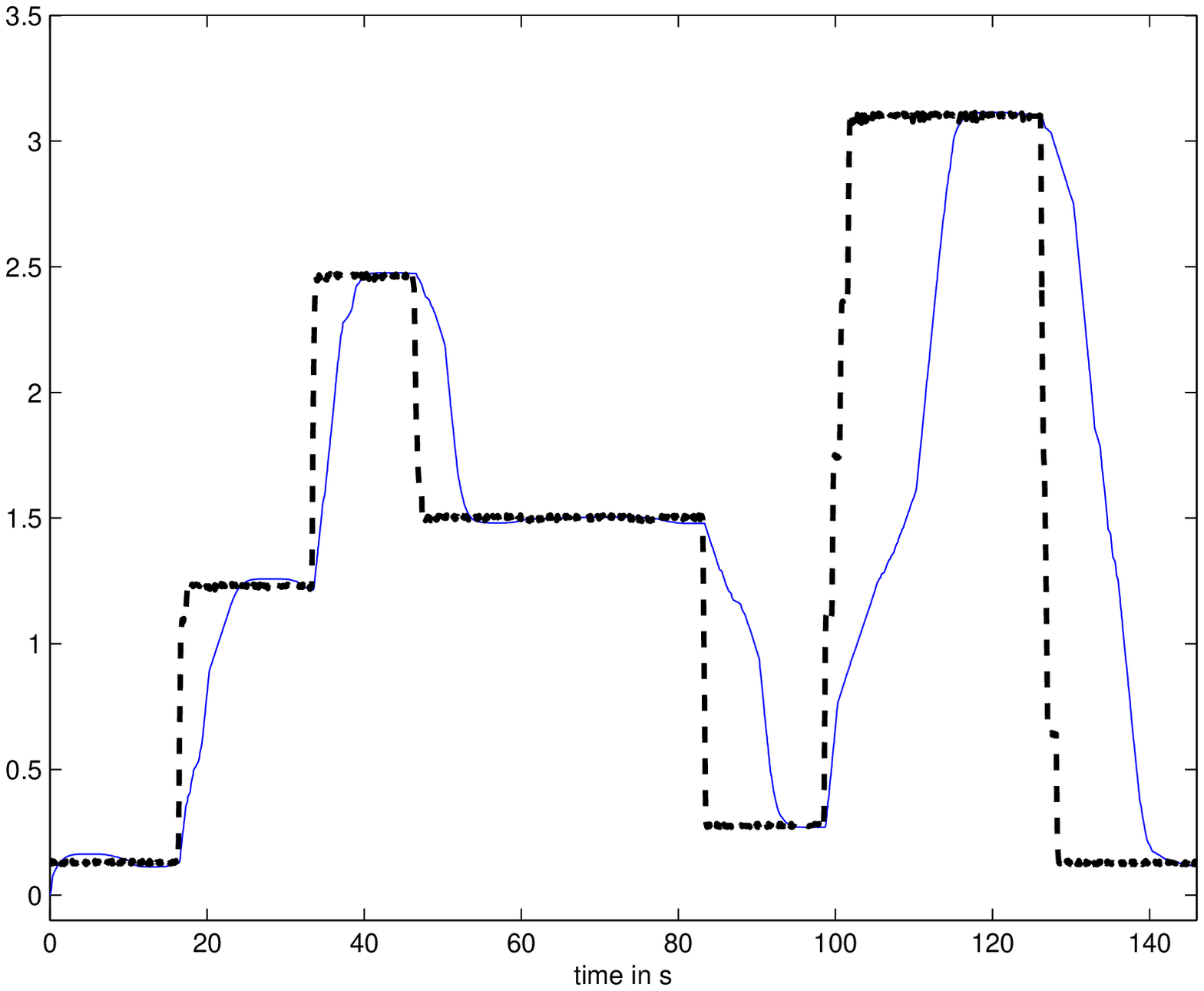}}}\\
\subfigure[Multiplicative power loss]{
\resizebox*{8cm}{!}{\includegraphics{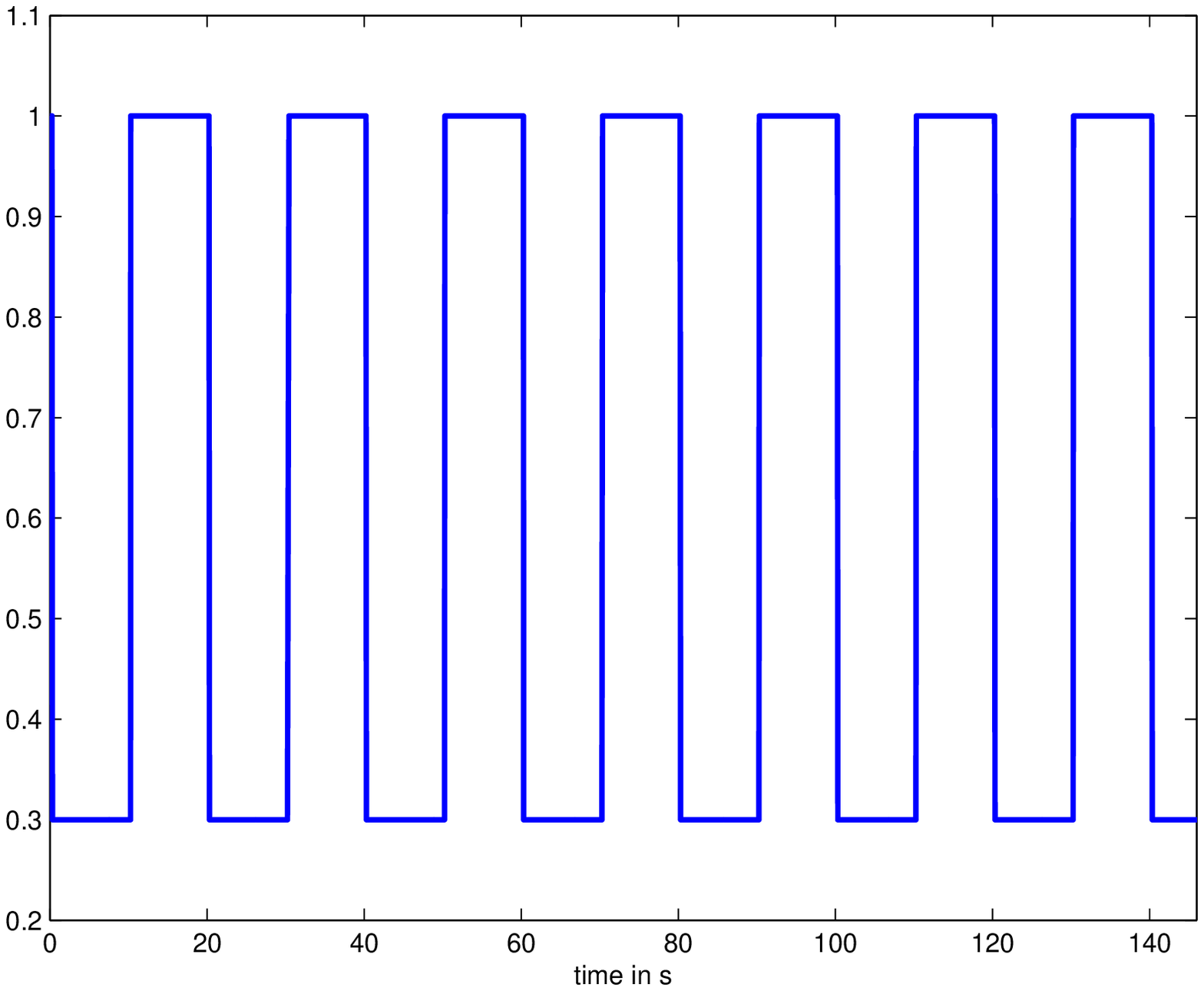}}}%
\caption{Experimental results: faulty case}%
\label{rez1ad}
\end{center}
\end{figure*}

Set
$$\bar v = \pi v$$ 
where 
\begin{itemize}
\item $\pi =1$ corresponds to the fault-free case,
\item Figure \ref{rez1ad}-(c) shows the case $0 \leq \pi \nleq 1$, which corresponds to a power loss of the actuator.
\end{itemize}
Figures \ref{rez1ad}-(a) and \ref{rez1ad}-(b) display an excellent fault tolerant control even with violent faults. 

\subsubsection{System 2}
\begin{figure*}[htb]
\begin{center}
\subfigure[Control (--,blue) and control limits (- -,red)]{
\resizebox*{8cm}{!}{\includegraphics{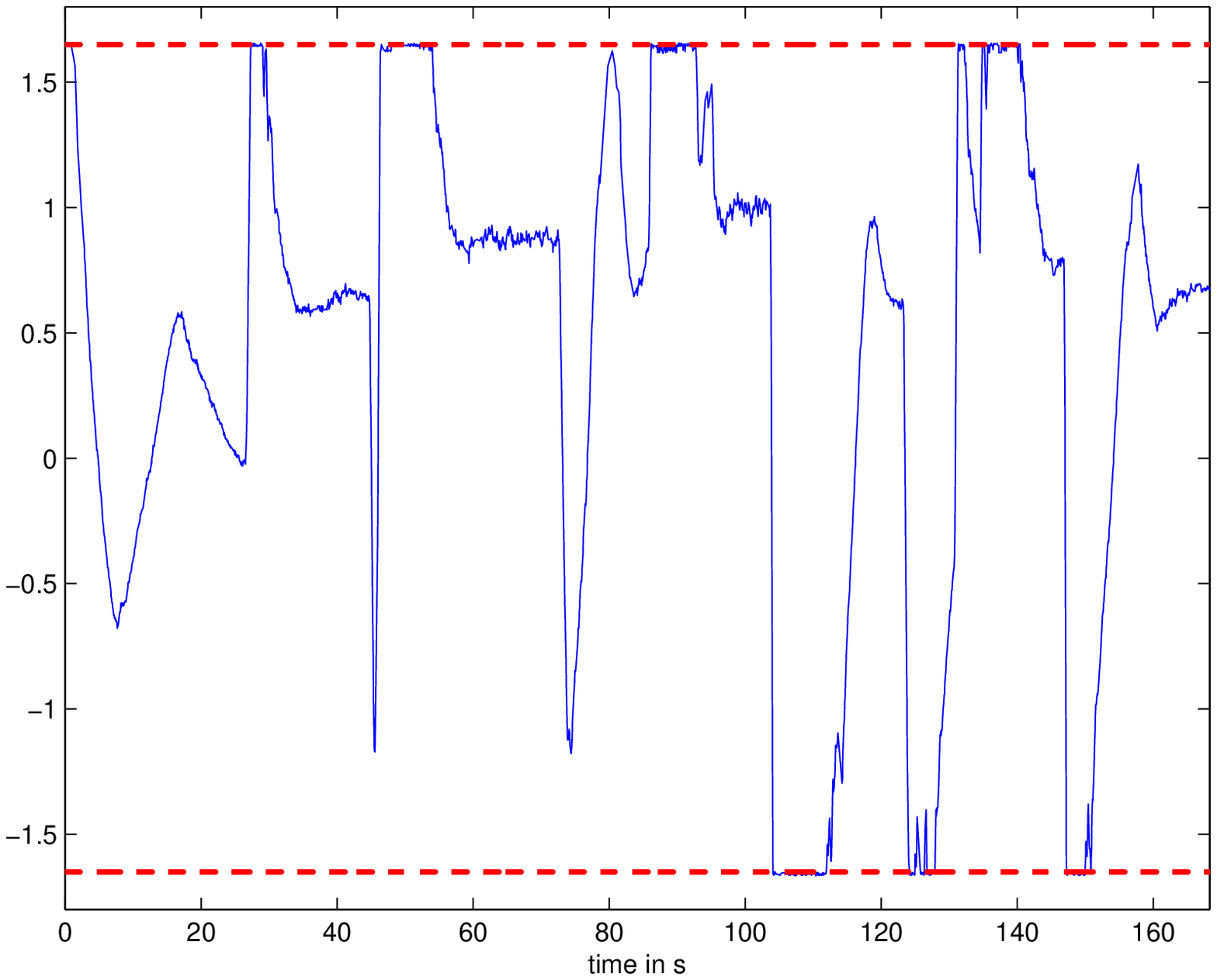}}}%
\subfigure[Setpoint (- - black) and Output (-- blue)]{
\resizebox*{8cm}{!}{\includegraphics{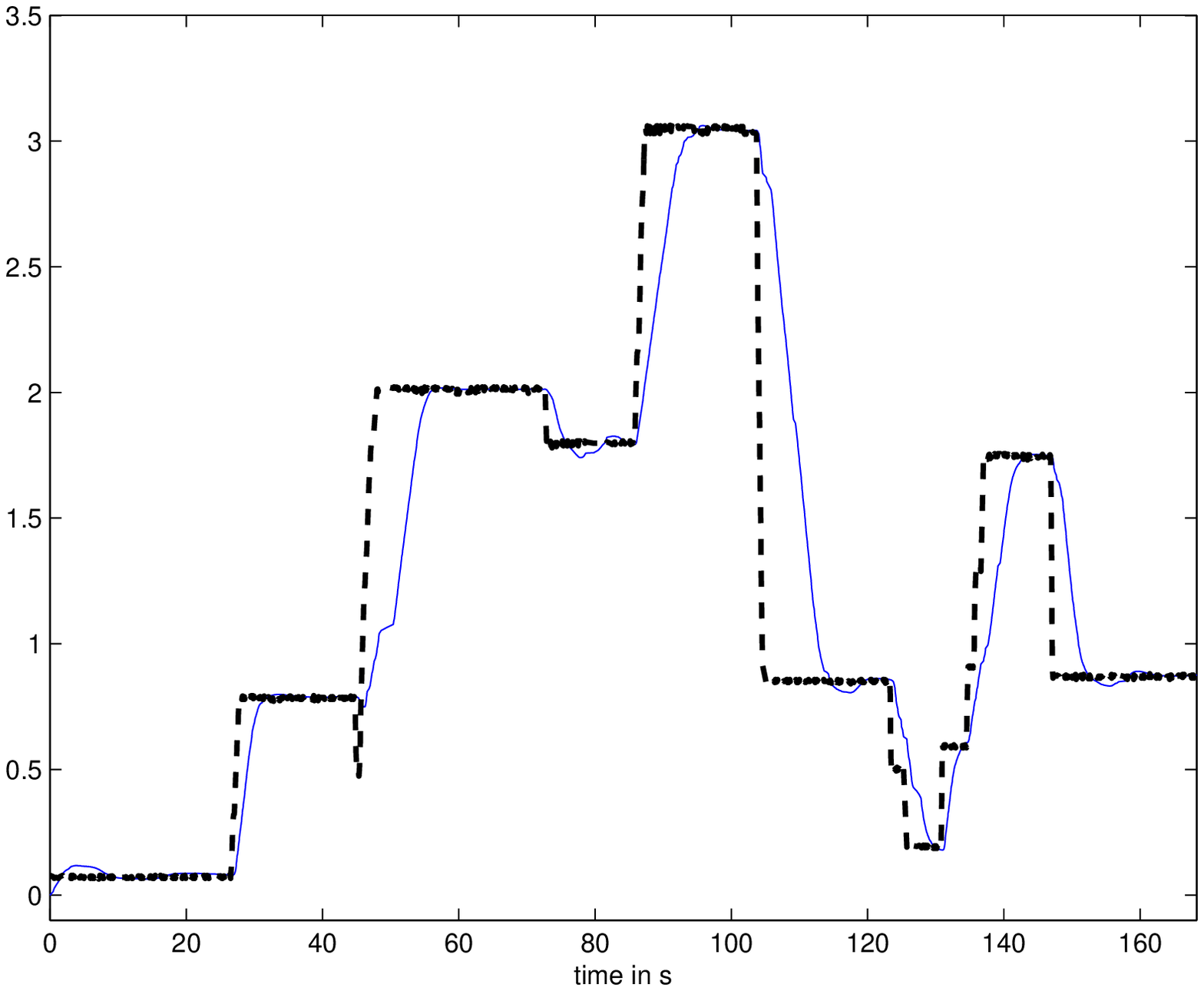}}}%
\caption{Experimental results: fault free case}%
\label{rez2}
\end{center}
\end{figure*}

Consider the nonlinear unstable system:
$$9\dot vv^2+v^3=\ddot y+\dot y +y$$
Close the loop with the same iP as above. Figure \ref{rez2}  displays excellent tracking and fault accommodation performances.

\begin{figure*}[htb]
\begin{center}
\subfigure[Control (-- blue) and control limits (- - red)]{
\resizebox*{8cm}{!}{\includegraphics{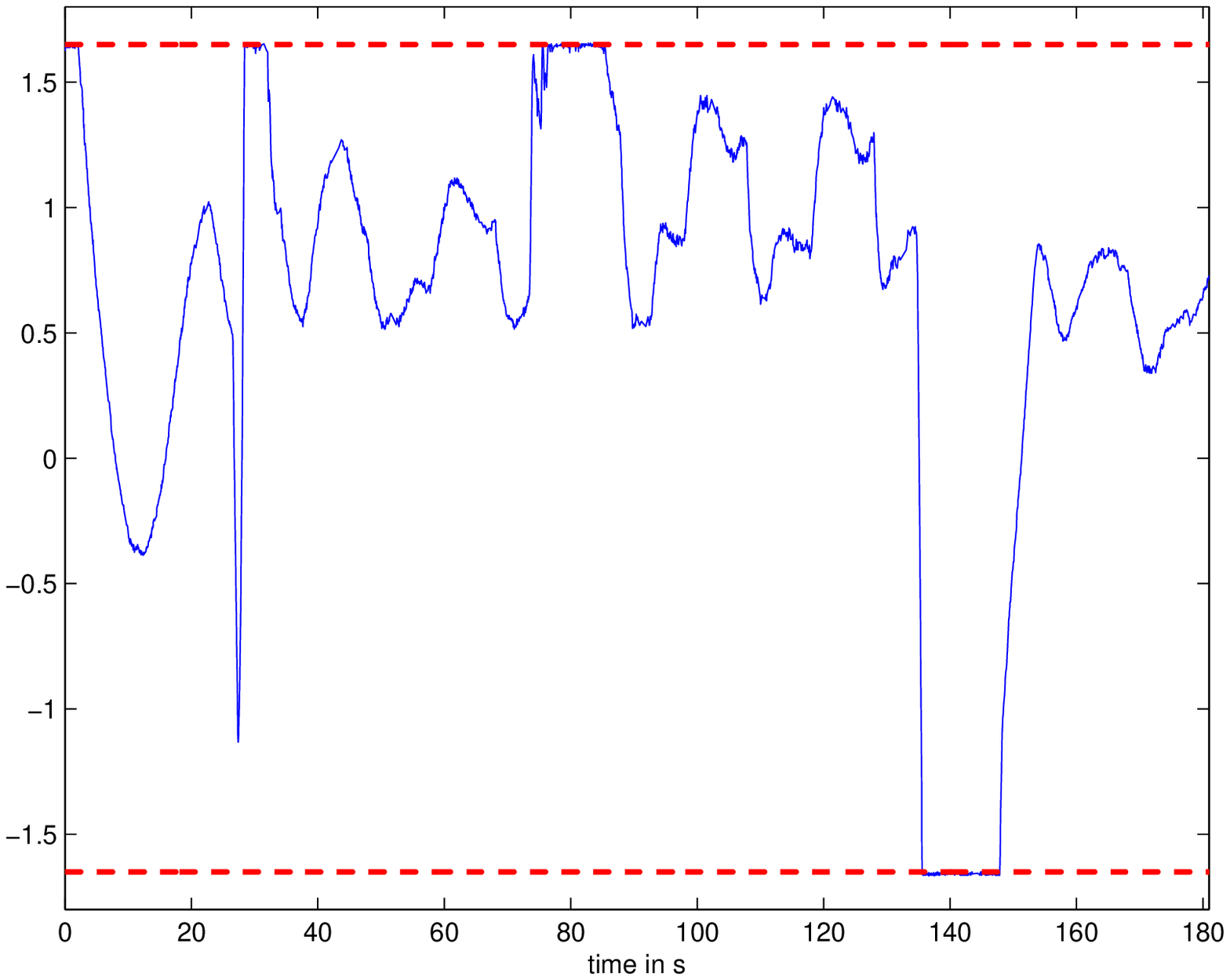}}}%
\subfigure[Setpoint (- - black) and output (-- blue)]{
\resizebox*{8cm}{!}{\includegraphics{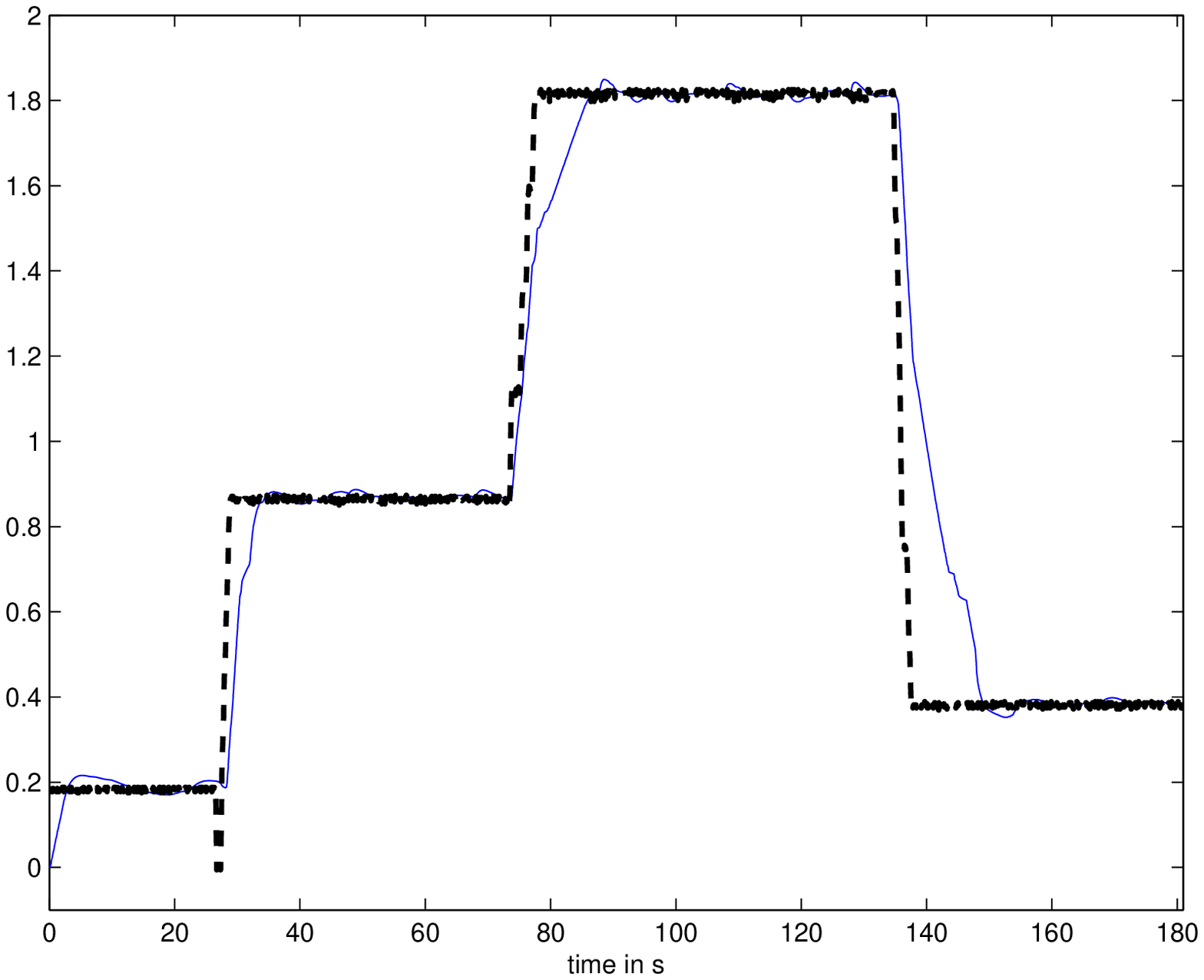}}}\\
\subfigure[Multiplicative power loss]{
\resizebox*{8cm}{!}{\includegraphics{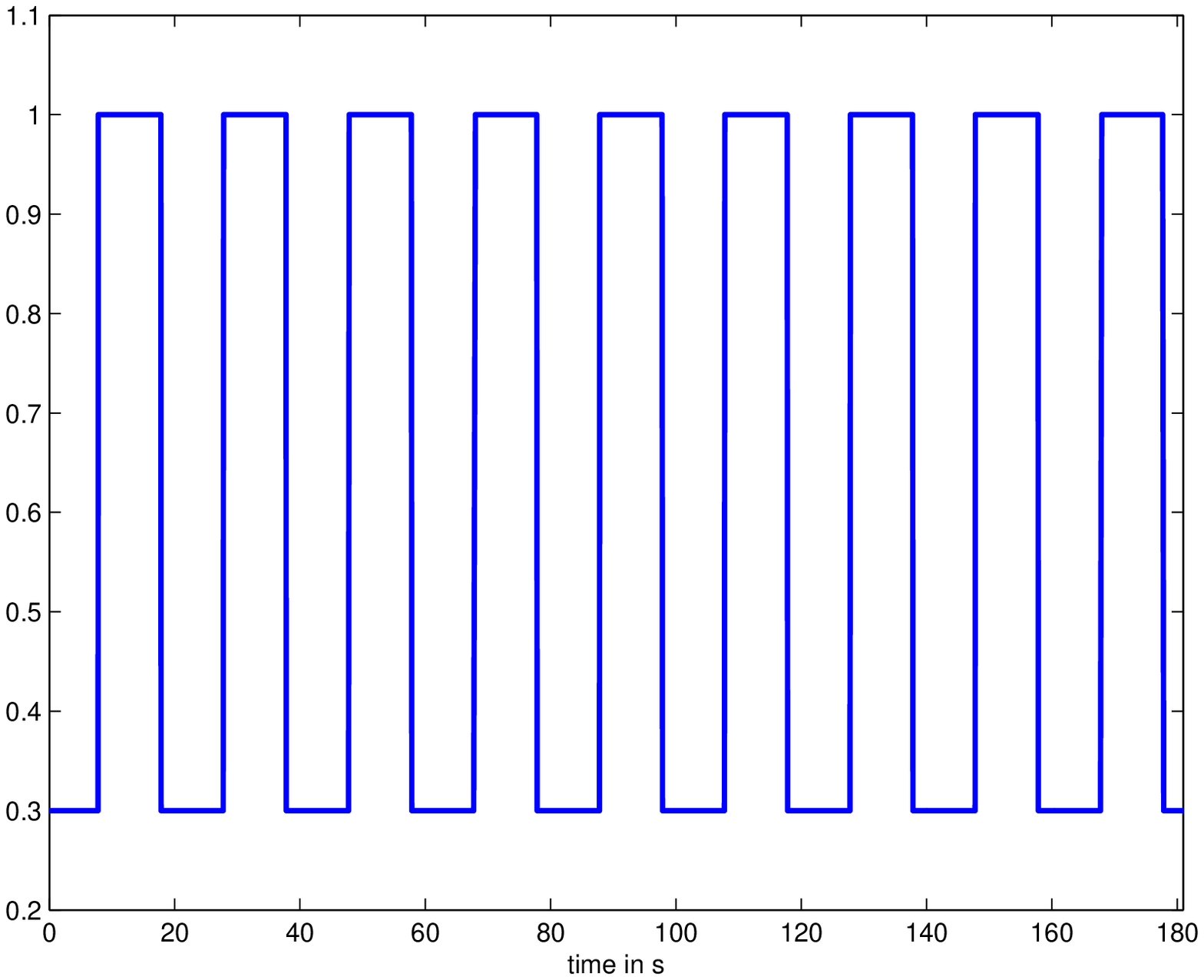}}}%
\caption{Experimental results: faulty case}%
\label{rez2ad}
\end{center}
\end{figure*}

\section{Conclusion}\label{conclu}
This communication has demonstrated that the intelligent controllers, which are associated to model-free control, may easily be implemented on cheap and small programmable devices. Future applications will show that 
our academic numerical experiments may be easily extended to more realistic case-studies. It should lead to great industrial opportunities for this new setting.

\vspace{0.3cm}
{\small 
\noindent {\bf Acknowledgement}: The authors thank the department of \textit{G\'enie \'Electrique et Informatique Industrielle},  \textit{Institut Universitaire de Technologie Nancy-Brabois},  \textit{Universit\'{e} de Lorraine}, 
for its most friendly help and its loan of some essential devices.}


\end{document}